
\input amstex
\documentstyle{amsppt} 


\def\cal{\Cal}
\define\U{\Cal U}
\define\M{\Cal M}
\define\V{\Cal V}
\define\E{\Cal E}
\define\F{\Cal F}
\define\W{\Cal W}
\redefine\O{\Cal O} 
\define\CB{\Cal B}  
\define\CC{\Cal C}  
\define\CQ{\Cal Q}  

\define\R{{\Bbb R}}
\define\Q{{\Bbb Q}}
\define\Z{{\Bbb Z}}
\define\C{{\Bbb C}}
\define\N{{\Bbb N}} 
\define\CP{{\Bbb P}}

\define\a{\alpha}
\undefine\b
\define\b{\beta}

\undefine\l
\define\l{\lambda}

\define\barPhi{{\bar\Phi}}

\define\Gc{G_{\Bbb C}}
\undefine\gg
\define\gg{\frak g}
\define\jgp{\frak p}  


\define\Sym{\operatornamewithlimits{Sym}}
\define\Ext{\operatornamewithlimits{Ext}}
\define\Hom{\operatornamewithlimits{Hom}}

\define\Res{\operatornamewithlimits{Res}} 
\define\apc{\operatorname{\hat{\ast}}}   
\define\apq{\operatorname{\hat{\cdot}}}  
\define\End{\operatorname{End}}          
\define\vspan{\operatorname{span}}       

\define\pd{\partial}
\define\rel{{}^r}
\undefine\Im
\define\Im{\text{\rom{Im}}}

\define\h{\hbar}         
\define\jdef#1{{\it #1}} 
\define\ti#1{\tilde{#1}} 
\define\jend#1{\End(#1)} 
\define\jpd #1{\frac{\partial}{\partial {#1}}} 
\define\qpsi{\gg^{\Psi}_-}


\TagsOnRight
\TagsAsMath
\loadmsbm

\topmatter

\abstract 
  In this article we describe three constructions of complex
  variations of Hodge structure, proving the existence of interesting
  opposite filtrations that generalize a construction of Deligne. We
  also analyze the relation between deformations of Frobenius modules
  and certain maximally degenerate variations of Hodge
  structures. Finally, under a certain generation hypothesis, we show
  how to construct a Frobenius manifold starting from a deformation
  of a Frobenius module.
\endabstract

\author Javier Fernandez \\ Gregory Pearlstein
\endauthor

\address Department of Mathematics, University of Utah, Salt Lake
City, Utah 84112-0090, USA.
\endaddress

\email jfernand\@math.utah.edu
\endemail

\address Department of Mathematics, University of California, Irvine,
California 92697-3875, USA.
\endaddress

\email gpearlst\@math.uci.edu
\endemail

\leftheadtext{Javier Fernandez and Gregory Pearlstein}

\title Opposite filtrations, variations of Hodge structure, \\
       and Frobenius modules
\endtitle

\rightheadtext{Opposite filtrations and Frobenius modules}

\dedicatory Dedicated to Professor Yuri I. Manin, on the occasion of
  his sixty fifth birthday 
\enddedicatory

\endtopmatter

\document



\head 1.\quad Introduction
\endhead

\par Let $X$ be a complex manifold.  Then, an \jdef{unpolarized
complex variation of Hodge structure $(E,\nabla,\F,\barPhi)$ over $X$}
consists of a flat, $C^\infty$ complex vector bundle $(E,\nabla)$ over $X$
equipped with a decreasing Hodge filtration $\F$ and an increasing
filtration $\barPhi$ such that
\roster
\item"(1.1)" $\F$ is holomorphic with respect to $\nabla$, and 
      $\nabla(\F^p)\subseteq \F^{p-1}\otimes\Omega^1_X$;
\item"(1.2)" $\barPhi$ is anti-holomorphic with respect to $\nabla$, and 
      $\nabla(\barPhi_q)\subseteq\barPhi_{q+1}\otimes\overline{\Omega^1_X}$;
\item"(1.3)"  $E = \F^p\oplus\barPhi_{p-1}$ for each index $p$ (i.e\. $\F$ is 
      opposite to $\barPhi$).
\endroster

\par Alternatively, conditions $(1.1)$--$(1.3)$ are equivalent to the 
assertion that the $C^{\infty}$ decomposition
$$
        E = \bigoplus_p\, \U^p                                    \tag{1.4}
$$
defined by the Hodge bundles $\U^p = \F^p\cap\barPhi_p$ satisfy Griffiths'
transversality:
$$
   \nabla:\E^0(\U^p) \to \E^{0,1}(\U^{p+1})\oplus\E^{0,1}(\U^p)\oplus
                      \E^{1,0}(\U^p)\oplus\E^{1,0}(\U^{p-1})        
                                                                \tag{1.5}
$$

\par In the present article, we describe three constructions of complex
variations of Hodge structure, and their relationship with quantum cohomology
and Frobenius manifolds.  More precisely, in \S 2 we recall that:
\roster
\item"(a)" Given a variation of graded-polarized mixed Hodge structure
$\V$, the Hodge filtration $\F$ of $\V$ pairs with the convolution
$$
      \barPhi = (\bar\F^{\vee}\ast\W)
$$
of the weight filtration $\W$ of $\V$ and $\bar\F^{\vee}_q =
\bar\F^{-q}$ to define an unpolarized $\C$VHS on the underlying
$C^{\infty}$ flat vector bundle of $\V$.

\item"(b)" Near a maximally unipotent boundary point, the $B$-model variation
of Hodge structure attached to a degenerating family of Calabi--Yau threefolds
gives rise to a variation of mixed Hodge structure of Hodge--Tate type, whose
extension data encodes the quantum cohomology of the mirror.

\item"(c)" The Higgs field $\theta$ obtained by application of construction
$\text{(a)}$ to $\text{(b)}$ also determines the quantum product of the 
mirror, i.e. in this case mirror symmetry can be understood as a duality of 
$\C$VHS.
\endroster

\par Following \cite{D3}, in \S 3 we generalize the construction
$(\text{b})$ by showing that given any admissible variation of
graded-polarized mixed Hodge structure $\V\to\Delta^{*n}$, there
exists a canonical filtration $\underline{\Psi}$ which is opposite to
the Hodge filtration $\F$ near zero, and hence defines an unpolarized
$\C$VHS.  Moreover, when the limiting mixed Hodge structure is
Hodge--Tate, this opposite filtration coincides with the relative
weight filtration of $\V$.

\vskip 3pt

\par The even cohomology of a compact smooth manifold $X$ of dimension $k$
is naturally a Frobenius algebra with respect to the cup product and
intersection pairing. For K\"ahler manifolds, the quantum product
provides a deformation of this algebra. These constructions restrict
to $V=\oplus_p H^{p,p}(X)$ if $X$ is Calabi-Yau. A Frobenius module of
weight $k$ is an abstract version of the module structure obtained
when the product operation on $V$ is restricted to a module over $\Sym
H^{1,1}(X)$. In \S 4 we review these notions as well as the
equivalence between families of Frobenius modules and variations of
Hodge structure with special degenerating behavior (Theorem~(4.14)).

\vskip 3pt

\par In general, the construction of a Frobenius module arising from an
algebra as described above implies a loss of information. Yet, in some
cases, it is possible to recover the full product structure. In \S 5
we discuss how this is the case when the weight of the module is at
most $5$ and, also, when the algebra is ``generated by $H^{1,1}$'', a
condition that has already been used in the context of quantum
cohomology~\cite{KM,Kr} as well as in the recent work of C. Hertling
and Yu. Manin~\cite{HM}.

\vskip 3pt

If a family of Frobenius modules generates a family of Frobenius
algebras, a natural question to ask is if it is possible to unfold this
last family into a Frobenius manifold. We study some cases where
families of Frobenius modules can be unfolded into (germs) of
Frobenius manifolds. We do this in two ways: an explicit construction
is given for low weight families (5.6), and a general argument
(Theorem~(5.8)) is presented using the techniques of Sections 2 and 3,
as well as~\cite{HM}.

\remark{Acknowledgment} 
  Part of the research presented in this paper was done at the Max
  Planck Institute and we are thankful for their hospitality and
  support. We would also like to thank Claus Hertling for very helpful
  discussions.
\endremark


\head 2.\quad Variations of Mixed Hodge Structure
\endhead

\par Let $V$ be a finite dimensional vector space.  Then, a decreasing
filtration of $V$ is an exhaustive sequence of subspaces
$$
         0\subseteq\cdots\subseteq F^p\subseteq F^{p-1}\subseteq\cdots
          \subseteq V
$$
of $V$ such that $F^p\subseteq F^{p-1}$ for all $p$.  Likewise, an 
increasing filtration $\Psi$ of $V$ is an exhaustive sequence of subspaces 
of $V$ such that $\Psi_q\subseteq \Psi_{q+1}$ for each index $q$.  As alluded
to in the introduction, a decreasing filtration $F$ is then said to be 
\jdef{opposite} to an increasing filtration $\Psi$ if and only if 
$$
      V = F^p\oplus\Psi_{p-1}                                   \tag{2.1}
$$
for each index $p$.  Similarly, given an integer $k$, a pair of decreasing
filtrations $F$ and $G$ are said to be $k$-opposed if and only if $F$ is
opposite to the increasing filtration $\Psi_q := G^{k-q}$.   

\definition{Definition 2.2} Let $A$ be a subfield of $\R$.  Then, a 
\jdef{pure $A$--Hodge structure of weight $k$} is pair $(V_A,F)$ consisting 
of a finite dimensional $A$-vector space $V_A$ and a decreasing Hodge 
filtration $F$ of $V_{\C} = V_A\otimes\C$ by complex subspaces such 
that the conjugate filtration $\bar F$ defined by the real structure 
$V_{\R} = V_A\otimes\R$ is \jdef{$k$-opposed} to $F$.
\enddefinition

\par By the Hodge decomposition theorem, the primitive\footnote{In
fact, the same is true without restricting to the primitive cohomology (see
Chapter 5, Section 6 of~\cite{We}).} $k$'th cohomology of a smooth
complex projective variety $X$ carries a pure Hodge structure of
weight $k$.  Accordingly, one defines a polarization of a pure Hodge
structure $(V_A,F)$ of weight $k$ to be a $(-1)^k$--symmetric bilinear
form $Q:V_A\otimes V_A\to A$ such that

\roster
\item $Q(F^p,F^{k-p+1}) = 0$;
\item $i^{p-q}Q(v,\bar v)$ is positive definite on $F^p\cap\bar F^{k-p}$;
\endroster
for each index $p$.

\definition{Definition 2.3} Let $S$ be a complex manifold and $A$ be a
subfield of $\R$.  Then, a \jdef{variation of pure, polarized
$A$--Hodge structure of weight $k$ over $S$} consists of a local
system of $\V_A$ of finite dimensional $A$-vector spaces over $S$
equipped with a decreasing Hodge filtration $\F$ of $\V =
\V_A\otimes{\cal O}_S$ by holomorphic subbundles, and a flat
$(-1)^k$-symmetric bilinear form $Q:\V_A\otimes\V_A\to A$ such that
\roster
\item $\F$ and $\bar\F$ are $k$-opposed;
\item $\F$ is horizontal, 
      i.e. $\nabla(\F^p)\subseteq \F^{p-1}\otimes\Omega^1_S$;
\item $Q$ polarizes each fiber of $\V$.
\endroster
\enddefinition

\par In particular, by the work of Griffiths~\cite{G}, each smooth projective
morphism $f:X\to S$ of complex algebraic varieties gives rise to a
variation of pure $\Q$--Hodge structure of weight $k$ on $\V_{\Q} =
R^k_{f*}(\Q)$. On the primitive part of $R^k_{f*}(\Q)$ the variation
is, also, polarized. Dropping the requirement that $f$ be smooth and/or
projective, one then obtains the notion of a variation of
graded-polarized mixed Hodge structure:

\definition{Definition 2.4~\cite{D2}} Let $A$ be a subfield of $\R$.
Then, an \jdef{$A$--mixed Hodge structure} consists of a finite
dimensional $A$-vector space equipped with a decreasing Hodge
filtration $F$ of $V_{\C}$ together with an increasing weight
filtration $W(V_A)$ of $V_A$ such that $F$ induces a pure $A$--Hodge
structure of weight $k$ on each non-trivial quotient $Gr^W_k =
W_k/W_{k-1}$ of the complexification $W$ of $W(V_A)$ via the rule
$$
        F^p Gr^W_k = \frac{F^p\cap W_k + W_{k-1}}{W_{k-1}}\qquad
$$ 
Likewise, a \jdef{graded-polarized mixed Hodge structure} is just a
mixed Hodge structure endowed with a choice of polarization $Q_k$ for
each non-trivial quotient $Gr^W_k$.
\enddefinition

\definition{Definition 2.5~\cite{SZ}} Let $S$ be a complex manifold
and $A$ be a subfield of $\R$. Then, a \jdef{variation of graded-polarized
$A$--mixed Hodge structure over $S$} consists of a local system $\V_A$
of finite dimensional $A$-vector spaces over $S$ equipped with a
decreasing Hodge filtration $\F$ of $\V = \V_A\otimes{\cal O}_S$ by
holomorphic subbundles, and an increasing weight filtration $\W(\V_A)$
of $\V_A$ by flat subbundles, together with a collection of flat,
non-degenerate bilinear forms $Q_k$ such that:
\roster
\item $\nabla(\F^p)\subseteq \F^{p-1}\otimes\Omega^1_S$;
\item $(Gr^{\W}_k(\V_A),\F Gr^{\W}_k,Q_k)$ is a variation of pure
      polarized $A$--Hodge structure of weight $k$.
\endroster
\enddefinition

\remark{Remark} In definitions $(2.2)$--$(2.5)$ one may replace $A$ by a
finitely generated noetherian subring of $\R$ such that $A\otimes\Q$ is a 
field.  However,  in this context, one only requires $W$ to be a filtration
of $V_A\otimes\Q$.
\endremark
\vskip 3pt

\par Let $\V$ be a variation of pure, polarized Hodge structure of weight 
$k$ and $\barPhi$ be the increasing filtration of $\V$ defined by the rule
$$
        \barPhi_p = \bar\F^{k-p}                                
$$
Then, the fact that $\bar\nabla = \nabla$ and $\F$ is holomorphic, horizontal 
and $k$-opposed to $\bar\F$ implies that $\barPhi$ is an anti-holomorphic
filtration of $\V$ which is opposite to $\F$ such that
$$
        \nabla(\barPhi_q)\subseteq\barPhi_{q+1}\otimes\bar\Omega^1
                                                                \tag{2.6}
$$
i.e. $(\F,\barPhi)$ defines an unpolarized $\C$VHS on the underlying flat 
$C^{\infty}$ bundle of $\V$.  Accordingly~\cite{Si}, one obtains an 
associated Higgs bundle structure $\bar\pd + \theta$ on the underlying
$C^{\infty}$ vector bundle of $\V$ as follows: Let
$$
        \V = \bigoplus_{p+q=k}\, {\cal H}^{p,q}
$$
denote the Hodge decomposition $(1.4)$ defined by the Hodge bundles
$$
        \U^p = {\cal H}^{p,q} = \F^p\cap \barPhi_p
$$
Then, by virtue of equation $(1.5)$, we can write
$$
        \nabla = \tau + \bar\pd + \pd + \theta                  \tag{2.7}
$$
where
$$
\aligned
        \tau&(\U^p)\subseteq\E^{0,1}(\U^{p+1}),                 \\
        \pd&(\E^0(\U^p))\subseteq\E^{1,0}(\U^p),
\endaligned\qquad
\aligned
        \bar\pd&(\E^0(\U^p))\subseteq\E^{0,1}(\U^p)             \\
        \theta&(\U^p)\subseteq\E^{1,0}(\U^{p-1})                
\endaligned                                                     \tag{2.8}
$$
Moreover, upon expanding out the flatness condition $\nabla^2 = 0$ relative
to $(2.7)$ and $(2.8)$, one obtains the Higgs field condition
$$
        (\bar\pd + \theta)^2 
        = \bar\pd^2 + \bar\pd\theta + \theta\wedge\theta = 0    \tag{2.9}
$$

\par More generally, mutatis mutandis, given an unpolarized $\C$VHS, one
obtains a corresponding Higgs bundle structure $\bar\pd + \theta$ on the
underlying $C^{\infty}$ vector bundle of the variation via the above
constructions.  
\vskip 3pt

\par Let $\V$ be a variation of graded-polarized mixed Hodge structure.
Then, as consequence of the following result, one obtains a system of 
Hodge bundles $(1.4)$ on the underlying $C^{\infty}$ bundle of $\V$
which satisfy the transversality condition $(1.5)$, and hence a 
corresponding Higgs bundle structure $\bar\pd + \theta$.

\proclaim{Theorem 2.10~\cite{D2,D4}} Let $(F,W)$ be a mixed Hodge structure.
Then, there exists a unique, functorial bigrading
$$
        V_{\C} = \bigoplus_{p,q}\, I^{p,q}                      \tag{2.11}
$$
of the underlying complex vector space $V_{\C}$ such that
\roster
\item"(a)" $F^p = \oplus_{a\geq p}\, I^{a,b}$ for each index $p$;
\item"(b)" $W_k = \bigoplus_{a+b\leq k}\, I^{a,b}$ for each index $k$;
\item"(c)" For each bi-index $(p,q)$,
$$
        \overline{I^{p,q}} = I^{q,p} \mod \bigoplus_{a<q,\,b<p} I^{a,b}
                                                                \tag{2.12}
$$
\endroster
Moreover,
$$
        \bigoplus_{r\leq p}\, I^{r,s} = \sum_k\, \bar F^{k-p}\cap W_k
                                                                \tag{2.13}
$$
\endproclaim
\demo{Proof} The existence of a direct sum decomposition $(2.11)$
satisfying conditions $(a)$--$(c)$ is shown in~\cite{D2}.  That this
decomposition also satisfies $(2.13)$ is shown in~\cite{D4}
(cf\.~\cite{P}).  Namely, in order to establish $(2.13)$, observe
that
$$
        \bigoplus_{r\leq p} I^{r,s} = \bigoplus_{r\leq p} \overline{I^{s,r}}
                                                                \tag{2.14}
$$
since, for any bi-index $(r,s)$ with $r\leq p$, the left hand side of $(2.14)$
contains every bi-index $(r',s') < (r,s)$.  Thus, it is sufficient to show
that 
$$
        \bigoplus_{r\leq p} \overline{I^{s,r}} 
        =  \sum_k\, \bar F^{k-p}\cap W_k                        \tag{2.15}
$$

\par By conditions $(a)$ and $(b)$, for each pair of indices $p$ and $k$,
$$
        \bar F^{k-p}\cap W_k = \overline{F^{k-p}\cap W_k}  
        = \bigoplus_{r+s\leq k,\, s\geq k-p}\, \overline{I^{s,r}}
                                                                \tag{2.16}
$$
and hence
$$
        \bar F^{k-p}\cap W_k \subseteq \bigoplus_{r\leq p} \overline{I^{s,r}}
                                                                \tag{2.17}
$$
since $r+s\leq k$ and $s\geq k-p \implies r+s\leq k\leq s+p\implies r\leq p$.
Thus, fixing $p$ and taking the sum of $(2.17)$ over all $k$ yields:
$$
        \sum_k\,\bar F^{k-p}\cap W_k 
        \subseteq \bigoplus_{r\leq p} \overline{I^{s,r}}        
$$
Conversely, by equation $(2.16)$,
$$
        \bigoplus_{r\leq p} \overline{I^{s,r}} 
        \subseteq \sum_k\,\bar F^{k-p}\cap W_k                  \tag{2.18}
$$
Indeed, given any bi-index $(s,r)$ with
$r\leq p$, there exists an index $k$ such that 
$
        r+s \leq k \leq s + p                           
$
(i.e. $r+s\leq k$ and $s\geq k-p$), and hence 
$\overline{I^{s,r}}\subseteq \bar F^{k-p}\cap W_k$.
\enddemo

\proclaim{Corollary 2.19} Let $\V$ be a variation of mixed Hodge structure,
and 
$$
        \V = \bigoplus_{p,q}\, {\Cal I}^{p,q}                   \tag{2.20}
$$
denotes the $C^{\infty}$ decomposition of $\V$ into a sum of $C^{\infty}$
subbundles defined by the pointwise application of Theorem $(2.10)$.
Then, the Hodge filtration $\F$ of $\V$ pairs with the increasing filtration
$$
        \barPhi_q = \sum_k\, \bar\F^{k-q}\cap\W_k               \tag{2.21}
$$
to define an unpolarized $\C$VHS for which the resulting $C^{\infty}$ 
decomposition $(1.4)$ is given by the formula
$$
        \V = \bigoplus_p\,\U^p,\qquad \U^p = \bigoplus_q\,{\Cal I}^{p,q}
                                                                \tag{2.22}
$$
\endproclaim
\demo{Proof} As in the pure case, the fact that $\F$ is holomorphic and 
horizontal, implies that $\barPhi$ is an anti-holomorphic filtration of $\V$
such that $\nabla(\barPhi_p)\subseteq\barPhi_{p+1}\otimes\bar\Omega^1$.  
Moreover, by Theorem $(2.10)$,
$$
        \F^p = \bigoplus_{a\geq p}\, {\Cal I}^{a,b},\qquad
        \barPhi_p = \bigoplus_{a\leq p}\, {\Cal I}^{a,b}        \tag{2.23}
$$
Consequently,
$
        \V = \F^p\oplus\barPhi_{p-1} = \bigoplus_{p,q}\, {\Cal I}^{p,q} 
$
and hence $(\F,\barPhi,\nabla)$ satisfy the axioms of an unpolarized
$\C$VHS.  Likewise, by equation $(2.23)$,
$
        \U^p := \F^p\cap\barPhi_p = \bigoplus_q\,{\Cal I}^{p,q}
$
\enddemo

\remark{Remark} Given a pair of increasing filtrations $A$ and $B$ of a 
vector space $V$, one defines the convolution $A\ast B$ to be the increasing
filtration of $V$ given by the rule
$$
        (A\ast B)_q = \sum_{r+s = q}\, A_r\cap B_s = \sum_k\, A_{q-k}\cap B_k
                                                                \tag{2.24}
$$
In particular, if for any decreasing filtration $F$ of $V$ one defines 
$F^{\vee}_r = F^{-r}$ then the increasing filtration $(2.21)$ is given
by the formula
$
        \barPhi = (\bar\F^{\vee}\ast\W)                         
$
\endremark
\vskip 3pt

\par Following~\cite{D3}, we now relate the above constructions with mirror
symmetry.  To this end, let $\V\to\Delta^{*r}$ be a variation of pure,
polarized Hodge structure of weight $k$ over a product of punctured
disks which has unipotent monodromy, and $T_j = e^{-N_j}$ denote
the monodromy of $\V$ about the $j$'th punctured disk.  Then, as a 
consequence of Schmid's orbit theorems~\cite{Sc}, one obtains an 
associated limiting mixed Hodge structure
$
        (F_{\infty},W[-k])                              
$
on the central fiber~\cite{D1} of the canonical extension of $\V$, where
$W[-k]_j = W_{j-k}$ denotes the shifted monodromy weight filtration
of $\V$:

\proclaim{Theorem 2.25~\cite{CK}} Let $\V\to\Delta^{*r}$ be a variation of
pure, polarized $A$--Hodge structure with unipotent monodromy, and 
$$
        {\Cal C} = \bigl\{\, \sum_j\, a_j N_j \mid a_j > 0 \,\bigr\}
$$
denote the corresponding monodromy cone.  Then, there exists a unique 
increasing filtration $\W(\V_A)$ of $\V_A$ by sub-local systems, called the 
monodromy weight filtration of $\V$, such that for every element $N\in\Cal C$:
\roster
\item $N(\W_j)\subseteq \W_{j-2}$; 
\item $N^j:Gr^{\W}_j @> \cong >> Gr^{\W}_{-j}$;
\endroster
\endproclaim

\par Suppose now that $\V\to\Delta^{*r}$ is a variation of pure Hodge
structure for which the corresponding limiting mixed Hodge structure
is Hodge--Tate, i.e. 
$$
        I^{p,q}_{(F_{\infty},W[-k])} = 0
$$ 
unless $p=q$.  Then, following~\cite{D3}, $\V$ should be regarded
as maximally degenerate.  Moreover, as consequence of the following
result, to each such maximally degenerate variation $\V$, one can
associate a corresponding Hodge--Tate variation $\V^{\circ}$:

\proclaim{Theorem 2.26~\cite{D3}} Let $\V\to\Delta^{*r}$ be a
variation of pure polarized Hodge structure of weight $k$ for which
the associated limiting mixed Hodge structure is Hodge--Tate.  Then,
the Hodge filtration $\F$ pairs with the shifted monodromy weight
filtration $\W[-k]$ of $\V$ to define a Hodge--Tate variation
$\V^{\circ}$ over a neighborhood of zero in $\Delta^{*r}$.
\endproclaim

\remark{Remark} By shrinking $\Delta^{*r}$ as necessary, we shall henceforth
assume that $(\F,\W)$ is Hodge--Tate over all of $\Delta^{*r}$.
\endremark
\vskip 3pt

\par To apply Theorem $(2.26)$ to the study of mirror symmetry, we now
recall the following definition:

\definition{Definition 2.27} Let $X$ be a polarized Calabi--Yau
manifold of dimension $d$, and $\overline{\M}_X$ be a smooth
partial compactification of the moduli space $\M_X$~\cite{B,Ti,To} of complex
structures on $X$, such that
$$
        D = \overline{\M}_X - \M_X         \tag{2.28}
$$ 
is a normal crossing divisor. Let $p\in D$ be such that, in a
neighborhood of $p$, $D=\sum_{j=1}^r D_j$ with $r= \dim \M_X = \dim
H^{d-1,1}(X)$ and $\{p\} = \cap_{j=1}^r D_j$. Then, $p$ is said to be
a \jdef{maximally unipotent boundary point of $\M_X$} if and only if
the variation of Hodge structure
$$
        \V_s = H^d(X_s)                                         \tag{2.29}
$$
over $\M_X$ has maximal unipotent monodromy on a neighborhood of $p$ in 
$\M_X$, i.e. the associated monodromy weight filtration $\W$ satisfies
the following two conditions:
\roster
\item $\dim\, Gr^{\W}_d = 1$, $\dim\, Gr^{\W}_{d-1} = 0$,  
      $\dim\, Gr^{\W}_{d-2} = r$;
\item $Gr^{\W}_{d-2} = \bigoplus_{j=1}^r\, N_j(Gr^{\W}_k)$.
\endroster
\enddefinition

\par Let $X$ be a smooth Calabi--Yau threefold, and $p$ be a maximally
unipotent boundary point of $\M_X$.  Then, it can be shown that the
limiting mixed Hodge structure of the variation $(2.29)$ at $p$ is
Hodge--Tate.  Moreover, according to~\cite{D3}, the resulting
variation $\V^{\circ}$ produced by Theorem $(2.26)$ should be viewed
as the $A$-model variation of the mirror $X^{\circ}$ of $X$.  

\vskip 3pt

\par To extract the quantum cohomology of $X^{\circ}$ using this approach, we 
recall~\cite{D3} that given a pair of variations of $\Z$--Hodge
structure $\Cal A$, and $\Cal B$ over $S$ of pure type $(p,p)$ and
$(p-1,p-1)$ respectively,
$$
        {\Ext}^1(\Cal A,\Cal B) = \Hom(\Cal A,\Cal B)\otimes{\Cal O}^*(S).
                                                                \tag{2.30}
$$
More generally, given a $\Z$-variation $\V$ of Hodge--Tate type, let $E_k$
denote the extension class $(2.30)$ attached to the extension
$$
        0\to Gr^{\W}_{2k-2}\to \W_{2k}/\W_{2k-4} \to Gr^{\W}_{2k}\to 0
$$

\proclaim{Lemma 2.31~\cite{D3}} Let $\V\to\Delta^{*r}$ be a
$\Z$-variation of Hodge--Tate type, and $\tilde\V$ denote the
canonical extension of $\V$ over $\Delta^r$ defined in~\cite{D1}.
Then, the Hodge filtration $\F$ of $\V$ extends holomorphically to
$\tilde\V$, and remains opposite to the weight filtration $\W$ over
$0\in\Delta^r$ if and only if
$$
        E_k\in \Hom(Gr^{\W}_{2k},Gr^{\W}_{2k-2})
               \otimes{\Cal O}^*_{mer}(\Delta^{*r})
$$
for each index $k$, where ${\Cal O}^*_{mer}(\Delta^{*r})$ denotes the group
of invertible holomorphic functions on $\Delta^{*r}$ which extend 
meromorphically to $\Delta^r$.
\endproclaim

\remark{Remark} The variation $\V^{\circ}$ produced by Theorem $(2.26)$
always satisfies the hypothesis of Lemma $(2.31)$.
\endremark
\vskip 3pt

\par Thus, given a smooth Calabi--Yau threefold $X$ and a maximal unipotent
boundary point $p$ of $\M_X$, one obtains via application of Lemma $(2.31)$ 
to the resulting Hodge--Tate variation $\V^{\circ}$, a corresponding set 
of extension classes $E_3$, $E_2$ and $E_1$ which encode the quantum product 
of $X^{\circ}$ as follows:  Let $\W$ be the weight filtration of $\V^{\circ}$.
Then, $\W$ is defined over $\Z$ with weight graded quotients:
$$
        Gr^{\W}_6(\V^{\circ}_{\Z}) \cong \Z,\qquad
        Gr^{\W}_4(\V^{\circ}_{\Z}) \cong \Z^r,\qquad
        Gr^{\W}_2(\V^{\circ}_{\Z}) \cong \Z^r,\qquad
        Gr^{\W}_0(\V^{\circ}_{\Z}) \cong \Z                     
$$
where $r = \dim H^{2,1}(X)$. Select a generator $1$ of
$Gr^{\W}_6(\V^{\circ}_{\Z})$.  Then, on account of maximal unipotent
monodromy, there exists a unique system of \lq\lq canonical
coordinates\rq\rq{} $(q_1,\dots,q_r)$ on the base of $\V^{\circ}$
relative to which the extension class $E_3$ assumes the form
$$
        E_3(1) = \sum_j\, q_j N_j(1)                            \tag{2.32}
$$
Moreover, as discussed in~\cite{D3}, upon expanding the logarithmic derivative
of $E_2$ relative to the system of canonical coordinates $(2.32)$, one 
obtains a generating function for the number of rational curves in the
mirror $X^{\circ}$.

\par Alternatively, as described in~\cite{P}, given a maximal
unipotent boundary point $p$ as above, both the canonical coordinates
$(2.32)$ and the quantum product of the mirror $X^{\circ}$ can be
described in terms of the Higgs field $\theta$ attached to
$\V^{\circ}$ by $(2.19)$ as follows: Let
$$
        H^{p,p} = Gr^{\W}_{6-2p}(\V^{\circ}_{\Z}),              
$$
$1$ be a generator of $H^{0,0}$, and $T_j = N_j(1)$.  Then, $\theta$
induces a map
$$
        H^{p,p}\to H^{p+1,p+1}\otimes\Omega^1(\Delta^{*r})      
$$
such that
\roster
\item"(a)"$\theta(\xi)1 = \sum_j\, T_j\otimes\Omega_j(\xi)$;
\item"(b)"$T_a\ast T_b = \theta(\xi_a)\circ\theta(\xi_b)\, 1$;
\endroster
where $(\Omega_1,\dots,\Omega_r)$ denotes the coframe of 
$\Omega^1(\Delta^{*r})$ defined by the 1-forms
$$
        \Omega_j = \frac{1}{2\pi i}\frac{dq_j}{q_j},
$$ 
$T_a\ast T_b$ represents~\cite{M,CoK} the quantum product on 
$H^{1,1}(X^{\circ})$, and $\Omega_j(\xi_k) = \delta_{jk}$.
\vskip 3pt

\par In \S 3, we shall extend the correspondence $\V\leftrightarrow\V^{\circ}$
given by Theorem $(2.26)$ to arbitrary variations of graded-polarized
mixed Hodge structures $\V\to\Delta^{*n}$ which are admissible in the
sense of~\cite{SZ} by constructing a suitable opposite filtration
$\Psi$ from the limiting mixed Hodge structure of $\V$.


\head 3.\quad Asymptotic Behavior
\endhead

\par Let $(E,\nabla)$ be a flat $\C$-vector bundle over $\Delta^{*n}$
with unipotent monodromy.  Then~\cite{D1}, up to isomorphism, there exists
a unique extension $\tilde E\to\Delta^n$ of $E$ such that $\nabla$
has at worst simple poles with nilpotent residues along the normal
crossing divisor
$$
        D = \Delta^n - \Delta^{*n} = \cup_j\, D_j               \tag{3.1}
$$

\par Equivalently, given a system of local coordinates $(s_1,\dots,s_n)$
on $\Delta^n$ relative to which the divisor $(3.1)$ assumes the form
$$
        s_1\cdots s_n = 0
$$
with $s_j = 0$ on $D_j$, and flat multivalued frame $(\sigma_1,\dots,\sigma_m)$
of $E$, the canonical extension $\tilde E$ described above can be identified 
with the locally free sheaf generated by the sections
$$
        \tilde\sigma_j 
        = e^{\frac{1}{2\pi i}\sum_j\, (\log s_j) N_j}\sigma_j
                                                                \tag{3.2}
$$
where $T_j = e^{-N_j}$ denotes the monodromy action of $E$ about $s_j=0$.
\vskip 3pt

\par In particular, if $\V\to\Delta^{*n}$ is a variation of graded-polarized
mixed Hodge structure for which the Hodge filtration $\F$ of $\V$ extends
holomorphically to the canonical extension $\tilde\V$, one can define a
corresponding nilpotent orbit $\F_{nilp}$ by simply extending $\F(0)$ to
a filtration of $\V$ which is constant with respect to the frame $(3.2)$.
\vskip 3pt

\par Alternatively, $\F_{nilp}$ may be described in terms of the 
associated period map
$$
        \varphi:\Delta^{*n}\to\M/\Gamma                         \tag{3.3}
$$      
defined by parallel translation of the Hodge filtration $\F$ of $\V$ to 
a fixed reference fiber $V = \V_{s_o}$ of $\V$ as follows:  
\vskip 3pt

\par Let $\V\to\Delta^{*n}$ be a variation of graded-polarized mixed Hodge
structure, and $V = \V_{s_o}$ denote a fixed reference fiber of $\V$. Let
$W$ denote the specialization of the weight filtration of $\V$ to $V$, 
and $Q = \{Q_k\}$ denote the corresponding specialization of the 
graded-polarizations of $\V$ to $Gr^W$.  Define $X$ to be the flag variety 
consisting of all decreasing filtrations $F$ of $V$ such that
$$
        \dim(F^p) = \text{rank}(\F^p)
$$
and let $\M$ be the classifying space~\cite{P} consisting of all filtrations 
$F\in X$ such that $(F,W)$ is a mixed Hodge structure which is 
graded-polarized by $Q$.  Then, the period map $(3.3)$ takes values in the 
quotient of $\M$ by the action of the monodromy group $\Gamma$ of $\V$.  

\proclaim{Theorem 3.4~\cite{P}} The classifying space $\M$
defined above is a complex manifold upon which the real Lie group
$$
        G = \{\, g\in GL(V)^W \mid Gr(g)\in Aut_{\R}(Q) \,\}    
$$
acts transitively by automorphisms, where $Gr(g):Gr^W\to Gr^W$ denotes
the map induced by an element of the stabilizer $GL(V)^W$ of $W$ on $Gr^W$.
\endproclaim

\par In particular, on account of Theorem $(3.4)$, the orbit
$$
        \check\M = \Gc.F_o \subseteq X                          \tag{3.5}
$$
of a point $F_o\in\M$ under the complex Lie group
$$
        \Gc = \{\, g\in GL(V)^W \mid Gr(g)\in Aut_{\C}(Q) \,\}  \tag{3.6}
$$
is well defined, independent of $F_o\in\M$.  

\remark{Warning} In general, $\Gc$ is not the complexification of $G$.
\endremark
\vskip 3pt

\par To continue, observe that the period map $(3.3)$ is locally liftable.
Accordingly, there exists a holomorphic map $F:U^n\to\M$ which makes the 
following diagram commute
$$
\CD
          U^n           @> F >>   \M                            \\
        @V p VV                  @VVV                           \\
      \Delta^{*n}  @> \varphi >> \M/\Gamma
\endCD                                                          \tag{3.7}
$$
where $U^n\subset\C^n$ denotes the product of upper half-planes on which
the imaginary parts of the standard coordinates $(z_1,\dots,z_n)$ are
positive, and $p:U^n\to\Delta^{*n}$ denotes the standard covering map 
defined by the coordinates $s_j = e^{2\pi i z_j}$.  
\vskip 3pt

\par Thus, on account of the commutativity of $(3.7)$, the map
$$
        \psi(z) = e^{-\sum_j\, z_j N_j}.F(z)                    \tag{3.8}
$$
satisfies the periodicity condition 
$$
        \psi(z_1,\dots,z_j+1,\dots,z_n) = \psi(z_1,\dots,z_n) 
$$
and hence descends to a well defined holomorphic map
$$
        \psi(s):\Delta^{*n}\to\check\M                          \tag{3.9}
$$

\proclaim{Lemma 3.10} Let $\V\to\Delta^{*n}$ be a variation of 
graded-polarized mixed Hodge structure.  Then, the Hodge filtration $\F$ of 
$\V$ extends holomorphically to the canonical extension $\tilde\V\to\Delta^n$
if and only if 
$$
        F_{\infty} := \lim_{s\to 0}\psi(s)                      \tag{3.11}
$$
exists.  
\endproclaim
\demo{Proof} Relative to the trivialization $(3.2)$, $\F$ coincides with
the filtration $(3.9)$.
\enddemo

\par Likewise, after unraveling the above definitions, on finds that the 
pull back $p^*(\F_{nilp})$ of $\F_{nilp}$ to $U^n$ coincides with the 
nilpotent orbit 
$$
        F_{nilp} = e^{\sum_j\, z_j N_j}.F_{\infty}              \tag{3.12}
$$

\remark{Remark} To be coordinate free, $(3.11)$ and $(3.12)$ should be viewed
as follows:  Let $(s_1,\dots,s_n)$ be a system of local coordinates on 
$\Delta^n$ which are compatible with the given divisor structure $(3.1)$,
and $(\l_1,\dots,\l_n)$ be the corresponding system of coordinates on 
$T_0(\Delta^n)$ defined by the basis vectors $e_j=(\frac{\pd}{\pd s_j})_0$.  
Then, the period 
map
$$
        \varphi_{nilp}(\l_1,\dots,\l_n) 
        = e^{\frac{1}{2\pi i}\sum_j\, \log(\l_j) N_j}.F_{\infty}
$$
determines a variation of mixed Hodge structure over the
complement of the divisor $\l_1\cdots\l_n = 0$, with monodromy action
$T_j = e^{-N_j}$ about $\l_j = 0$, which is well defined, independent 
of the choice of local coordinates $(s_1,\dots,s_n)$ as above.
\endremark
\vskip 3pt

\par For variations of graded-polarized mixed Hodge structure, the analog
of the monodromy weight filtration $(2.25)$ is the relative weight
filtration $\rel W = \rel W(N,W)$ discussed in~\cite{SZ}.  Moreover,
based upon the study of degenerating families of varieties, Steenbrink
and Zucker proposed the following, now standard, definition of an
admissible variation of graded-polarized mixed Hodge structure over
$\Delta^*$:

\definition{Definition 3.13~\cite{SZ}} Let $\V\to\Delta^*$ be a variation of 
graded-polarized mixed Hodge structure with unipotent monodromy. Then, 
$\V$ is \jdef{admissible} if 
\roster
\item"(a)" The limiting Hodge filtration $(3.11)$ exists;
\item"(b)" The relative weight filtration $\rel W = \rel W(N,W)$ exists.
\endroster
\enddefinition

\par The admissibility conditions $(3.11)$ always hold in the pure case
as a consequence of Schmid's nilpotent orbit theorem.  For
multivariable degenerations, one defines admissibility via curve test
using $(3.13)$~\cite{K}.  Moreover, one has the following result:

\proclaim{Theorem 3.14~\cite{K}} Let $\V\to\Delta^{*n}$ be an admissible
variation of graded-polarized mixed Hodge structure with unipotent monodromy,
and 
$$
        {\Cal C} = \bigl\{\, \sum_j\, a_j N_j \mid a_j > 0 \,\bigr\}
$$
denote the corresponding monodromy cone.  Let $W$ denote the specialization
of the weight filtration of $\V$ to some fixed reference fiber $V = \V_{s_o}$.
Then,
\roster
\item $\rel W(N,W)$ exists for every element $N\in\Cal C$;
\item $\rel W = \rel W(N,W)$ is well defined, independent of $N\in\Cal C$;
\item $(F_{\infty},\rel W)$ is a mixed Hodge structure;
\item $N_1,\dots, N_n$ are $(-1,-1)$-morphisms of $(F_{\infty},\rel W)$.
\endroster
\endproclaim

\par Let $\V\to\Delta^{*n}$ be an admissible variation of graded-polarized
mixed Hodge structure with unipotent monodromy, and 
$$
        V = \bigoplus_{p,q}\, I^{p,q}                           \tag{3.15}
$$
denote the corresponding decomposition $(2.11)$ defined by the limiting
Hodge filtration $F_{\infty}$ and the relative weight filtration $\rel W$
of $\V$.  Define,
$$
        \Psi_p = \bigoplus_{a\leq p}\, I^{a,b}                  \tag{3.16}
$$
and
$
        \qpsi = \{\, \a\in\gg_{\C} 
                        \mid \a(\Psi_p)\subseteq \Psi_{p-1} \,\}
$
denote the subalgebra of $\gg_{\C} = Lie(\Gc)$ consisting of those elements
which preserve the increasing filtration $\Psi$ defined by $(3.16)$, and
act trivially on each layer $Gr^{\Psi}_p = \Psi_p/\Psi_{p-1}$ of 
$Gr^{\Psi}$.

\proclaim{Lemma 3.17} $\Psi$ is opposite to $F_{\infty}$. Moreover,
relative to the decomposition
$$
        \gg_{\C} = \bigoplus_{r,s}\, \gg^{r,s}                  \tag{3.18}
$$
induced by the bigrading $(3.15)$, 
$$
        \qpsi = \bigoplus_{r<0}\, \gg^{r,s}                  \tag{3.19}
$$
\endproclaim
\demo{Proof} That $\Psi$ is opposite to $F_{\infty}$ is a simple consequence
of definition $(3.16)$, and the fact [cf\. Theorem $(2.10)$] that 
$F_{\infty}^p = \oplus_{a\geq p}\, I^{a,b}$.  Likewise, since 
$$
        \gg^{r,s}(I^{a,b}) \subseteq I^{a+r,b+s}
$$
the sum $\oplus_{r<0}\, \gg^{r,s}$ maps $\Psi_p$ to $\Psi_{p-1}$, and
hence is contained in $\qpsi$.  Conversely, by equation $(3.16)$, if
$\a\in \qpsi$ then the components $\a^{r,s}$ with respect to the
decomposition $(3.18)$ can be non-zero only if $r<0$.
\enddemo

\proclaim{Corollary 3.20} Let $\V\to\Delta^{*n}$ be an admissible variation
with unipotent monodromy.  Then, on a neighborhood of the origin, the 
associated function $(3.9)$ admits a unique representation of the form
$$
        \psi(s) = e^{\Gamma(s)}.F_{\infty}                      \tag{3.21}
$$
with respect to $\qpsi$--valued holomorphic function $\Gamma(s)$
which vanishes at the origin.
\endproclaim
\demo{Proof} By $(3.19)$, $\qpsi$ is a vector space complement to
$
        Lie(\Gc^{F_{\infty}})
$ 
in $\gg_{\C}$.  Consequently, the map $u\mapsto e^u.F_{\infty}$ is a 
biholomorphism from a neighborhood of zero in $\qpsi$ onto a neighborhood
of $F_{\infty}$ in $\check\M$.   Accordingly, $\psi(s)$ admits a unique
representation of the type described above on a neighborhood of zero in
$\Delta^n$.
\enddemo

\par A priori, the filtration $\Psi$ defined above depends on the choice
of coordinates used in the construction of the limiting Hodge filtration
$(3.11)$.  However, as the following result shows, this is in fact not 
the case, since the monodromy logarithms $N_1,\dots,N_n$ preserve $\Psi$.

\proclaim{Theorem 3.22} $\Psi$ is independent of the choice of 
coordinates used in the definition of $F_{\infty}$.  Moreover,
$$
        \Psi = (\overline{F_{nilp}^{\vee}})\ast(\rel W)
             = (\overline{F_{\infty}^{\vee}})\ast(\rel W)       \tag{3.23}
$$
\endproclaim
\demo{Proof} Without loss of generality, any two systems of local coordinates
$(s_1,\dots,s_n)$ and $(\tilde s_1,\dots,\tilde s_n)$ on $\Delta^n$ which are 
compatible with the divisor structure $(3.1)$ may be  assumed to be of the 
form 
$
        \tilde s_j = f_j s_j
$
for some collection of holomorphic functions $f_1,\dots, f_n$ which do not
vanish at the origin.  Moreover, direct calculation shows that, under 
such changes of coordinates, the corresponding filtrations $(3.11)$ are
related by the equation:
$$
        \tilde F_{\infty} 
        = e^{-\frac{1}{2\pi i}\sum_j\, \log(f_j(0))N_j}.F_{\infty}
                                                                \tag{3.24}
$$

\par Let
$
        {\Cal N} = \text{span}_{\C}(N_1,\dots,N_n)
$.
Then,
$
        {\Cal N}\subseteq \gg^{-1,-1}                   
$
because each $N_j$ is a $(-1,-1)$-morphism of $(F_{\infty},\rel W)$.
Accordingly, 
$$
        N\in{\cal N} \implies 
        I^{a,b}_{(e^N.F_{\infty},\rel W)} 
        = e^N.I^{a,b}_{(F_{\infty},\rel W)}                     \tag{3.25}
$$
and hence
$$
        \Psi_{(\tilde F_{\infty},\rel W)} 
        = e^{\tilde N}.\Psi_{(F_{\infty},\rel W)}
$$
where 
$$
        \tilde N = -\frac{1}{2\pi i}\sum_j\, \log(f_j(0))N_j     \tag{3.26}
$$
On the other hand, 
$$
        e^{\tilde N}.\Psi_{(F_{\infty},\rel W)} 
        = \Psi_{(F_{\infty},\rel W)}                            
$$
since $\cal N\subseteq\gg^{-1,-1}$ and $\gg^{-1,-1}\subseteq \qpsi$ by 
Lemma $(3.17)$.  Consequently, 
$$
        \Psi = \Psi_{(F_{\infty},\rel W)}
             = e^{\tilde N}.\Psi_{(F_{\infty},\rel W)}
             = \Psi_{(\tilde F_{\infty},\rel W)}                \tag{3.27}
$$
is independent of the choice of coordinates used in the construction of 
the limiting Hodge filtration $(3.11)$.
\vskip 3pt
  
\par To verify  $(3.23)$, observe that 
$
        \Psi = (\overline{F_{\infty}^{\vee}})\ast(\rel W)
$ 
by virtue of equations $(3.16)$ and $(2.13)$.  To show that 
$$
        \Psi = (\overline{F_{nilp}^{\vee}})\ast(\rel W)         
$$ 
observe that $\cal N$ is closed under complex conjugation since $\bar
N_j = N_j$ for all $j$.  Consequently, when evaluated at any
particular point in $U^n$, $\overline{F_{nilp}} = e^N.\bar F_{\infty}$
for some element $N\in\cal N$.  On the other hand, by
definition~\cite{SZ}, $\cal N$ preserves the relative weight
filtration $\rel W$.  Thus,
$$
        (\overline{F_{nilp}^{\vee}})\ast(\rel W)
        = (e^N.\bar F_{\infty}^{\vee})\ast(\rel W)
        = e^N.(\bar F_{\infty}^{\vee})\ast(\rel W) 
        = e^N.\Psi = \Psi
$$
\enddemo

\par Finally, as a consequence of the above results, we obtain the following
generalization of Theorem $(2.26)$:

\proclaim{Theorem 3.28} 
  Let $\V\to\Delta^{*n}$ be an admissible variation of
  graded-polarized mixed Hodge structure with unipotent monodromy.
  Then, $\Psi$ extends to a filtration $\underline{\Psi}$ of $\V$ by
  flat subbundles which pairs with the Hodge filtration $\F$ of $\V$
  to define an unpolarized $\C$VHS on a neighborhood of the origin.
\endproclaim
\demo{Proof} Since each $N_j$ preserves $\Psi$, and $\V$ has only local 
monodromy, $\Psi$ extends to a filtration $\underline{\Psi}$ of $\V$ by 
flat subbundles.  To see that $\underline{\Psi}$ is opposite to $\F$, 
observe that after pulling everything back to the upper half-plane, we can 
write
$$
        F(z) = e^{\sum_j\, z_j N_j}\psi(s) 
             = e^{\sum_j\, z_j N_j}e^{\Gamma(s)}.F_{\infty}     \tag{3.29}
$$
relative to a $\qpsi$-valued holomorphic function $\Gamma(s)$ by
virtue of equations $(3.8)$ and $(3.21)$.  Moreover, since $\qpsi$ 
is nilpotent, we can write
$$
        e^{\sum_j\, z_j N_j}e^{\Gamma(s)} = e^{X(z)}            \tag{3.30}
$$  
for some function $X(z)$ with values in $\qpsi$.  Thus,
$$
\aligned
        V &= F^p_{\infty}\oplus\Psi_{p-1} 
           = e^{X(z)}.(F^p_{\infty}\oplus\Psi_{p-1})                    \\
          &= (e^{X(z)}.F^p_{\infty})\oplus(e^{X(z)}.\Psi_{p-1})
           = F^p(z)\oplus\Psi_{p-1}
\endaligned
$$
since $F(z) = e^{X(z)}.F_{\infty}$ and $\qpsi$ preserves $\Psi$.
\enddemo


\head 4.\quad Frobenius Modules
\endhead

Let $X$ be a compact K\"ahler manifold. Then the even part of
$H^{*}(X,\C)$ is naturally a Frobenius algebra with product given by
the cup product and the bilinear pairing arising from the intersection
form\footnote{The full $H^{*}(X,\C)$ is a $\Z_2$-graded Frobenius
algebra, as considered by Kontsevich and Manin in~\cite{KM}.}. The
quantum cohomology defines a deformation of this (``classical'')
structure parameterized by the complexified K\"ahler cone of $X$. In
this case, the deformed (``quantum'') product is defined in terms of a
function known as the Gromov-Witten potential.

In the case where $X$ is also Calabi-Yau, the cohomology decomposes
$H^*(X,\C) = \oplus_{p,q} H^{p,q}$ and the quantum product preserves
this bigrading. In particular, $\oplus_p H^{p,p}$ inherits a Frobenius
algebra structure. We should mention that these spaces also satisfy
nondegeneracy and positivity conditions encoded in the Hard Lefschetz
Theorem and the Hodge-Riemann relations.

\par In his analysis of mirror symmetry, D. Morrison constructed a
polarized variation of pure Hodge structure using the data described
above. This variation is known as the \jdef{$A$-model variation}
(\cite{M}, \cite{CoK,Chapter 8}). A close inspection of his
construction shows that not all of the algebra structure is used: it
is only the $\Sym H^{1,1}$-module structure on $\oplus_p H^{p,p}$ that
is needed. This observation is the starting point for the construction
of an equivalence between this type of structure, called a Frobenius
module, and certain maximally degenerating variations of pure Hodge
structures. In what follows we will describe this
construction. See~\cite{CF2} for further details.

For $k\in \N$ let $V=\oplus_{p=0}^k V_{2p}$ be a graded finite
dimensional $\C$-vector space and $\CB$ a symmetric non-degenerate
bilinear form on $V$ pairing $V_{2p}$ with $V_{2(k-p)}$. Let
$\{T_a\}_{0\leq a\leq m}$ be a $\CB$-self dual, graded basis of $V$.
We will refer to $\{T_a\}$ as an \jdef{adapted basis}.  For $0\leq
a\leq m$ define $\delta(a)$ by $\CB(T_{\delta(a)},T_b) = \delta_{a b}$
---the right hand side $\delta$ is Kronecker's $\delta$---
for all $b=0,\ldots,m$.  We also set $\ti{a} := p$ if and only if $
T_a \in V_p$ and assume that the map $\sim\ : \{0,\ldots,m\}
\rightarrow \{0,\ldots, 2k\}$ is increasing.

\definition{Definition 4.1} $(V,\CB,e,*)$ is a \jdef{graded
$V_2$-Frobenius module} of weight $k$ if
\roster     
\item $e\neq 0$ and $V_0 = \langle e \rangle$.
\item $V$ is a graded $\Sym V_2$-module under $*$.
\item For all $v_1, v_2 \in V$ and $w \in V_2$
$$
      \CB(w * v_1, v_2) \ =\ \CB(v_1,w * v_2)    \tag{4.2}
$$
\item $w * e = w$ for all $w\in V_2$.
\endroster
\enddefinition

Since $T_0\in V_0$, it must be a non-zero multiple of $e$ and we
assume that an adapted basis satisfies $T_0 = e$.  Clearly, the fact
that $V$ is a $\Sym V_2$-module is equivalent to
$$
  T_j * (T_l * T) \ =\ T_l * (T_j * T) \text{ for all }
  T_j, T_l \in V_2 \text{ and } T\in V.                 \tag{4.3}
$$

We say that $V$ is \jdef{real} if $V$ has a real structure, $V_\R$,
compatible with its grading, $*$ is real, $e\in V_\R$, and $\CB$ is
defined over $\R$.

To any real Frobenius module we can associate a Hodge--Tate mixed Hodge
structure, split over $\R$, whose canonical bigrading is  

$$
  I^{p,p} \ :=\  V_{2(k-p)}.                   \tag{4.4}
$$

The multiplication operator $L_w\in\jend{V}$, $w\in V_2$, is an
infinitesimal automorphism of the bilinear form
$$
  Q(v_a,v_b) \ :=\ (-1)^{k+\ti{a}/2} \CB(v_a,v_b),      \tag{4.5}
$$
as well as a $(-1,-1)$-morphism of the associated mixed Hodge
structure.  We will say that $w\in V_2\cap V_\R$ \jdef{polarizes} $V$
if the mixed Hodge structure $(I^{*,*},Q,L_w)$ is polarized~\cite{CK}.  
A real Frobenius module $V$ is said to be \jdef{polarizable} if it contains a
polarizing element. Given a polarizing element $w$, the set of
polarizing elements is an open cone in $V_{2} \cap V_\R$. We can then
choose a basis $T_1, \ldots, T_r$ of $V_{2}\cap V_\R$ spanning a
simplicial cone $\CC$ contained in the closure of the polarizing cone
and with $w\in \CC$.  Such a choice of a basis of $V_2$ will be called
a \jdef{framing} of the polarized Frobenius module.

\example{Example 4.6}
  If $X$ is a compact K\"ahler manifold of dimension $k$, let
  $V_{2p}:= H^{p,p}(X)$, $\CB_{int}$ the intersection pairing on
  $V:=\oplus_{p=0}^k V_{2p}$, and $\smallsmile$ the restriction of the cup
  product to $V$. Then, $(V,\CB_{int},1,\smallsmile)$ defines a polarizable
  Frobenius module. The real structure is induced by $H^*(X,\R)$.
\endexample

Given an adapted basis $\{T_0,\ldots,T_m\}$ of $V$, let
$z_0,\ldots,z_m$ be the corresponding linear coordinates on $V$ and
set $q_j:= \exp(2\pi i z_j)$ for $j=1,\ldots, r:=\dim V_2$ . If $U$ is
the upper-half plane, we may identify $U^r \cong (V_2 \cap V_\R)
\oplus i\,\CC$ and view the correspondence
$$
  \sum_{j=1}^r z_j T_j \in (V_2 \cap V_\R) \oplus i\,\CC
  \mapsto (q_1,\dots,q_r) \in (\Delta^*)^r.
$$

Let $V$ be a framed Frobenius module of weight $k$. Then the action of
$V_2$ on $V$ can be recovered from a homogeneous polynomial $\phi_0
\in \C[z_0,\ldots,z_m]$ of degree three, called
the \jdef{classical potential}. Indeed, if we let 
$$
  \phi_0(z_0,\ldots,z_m) \ :=\ \sum_{\ti{j}=2,\ 0\leq\ti{a},\ti{b}\leq
  2k} z_j z_a z_b\,\frac{C(\ti{a})}{12}\,\CB(T_j * T_a, T_b)\,, 
$$
with
$$
  C(\ti{a})\ :=\ \left\{ 
    \aligned
    & 2 \text{ if } k=3 \text{ and } \ti{a}=2,\\
    & 3 \text{ if } k\neq 3 \text{ and } \ti{a}=2 \text{ or }
    \ti{a}=2k-4,\\
    & 6 \text{ otherwise},
    \endaligned
  \right.
$$
then we recover the $V_2$-action by:
$$
  T_j * T_a \ :=\ \sum_{\ti{c}=\ti{a}+2}
  \frac{\partial ^3 \phi_0}{\partial z_j \partial z_a \partial
  z_{\delta(c)}} T_c\ ;\quad j=1,\dots,r\,.
$$

\example{Example 4.7}
  In the case of weight $k=3$, the classical potential is
$$
    \phi_0(z_0,\ldots,z_m) = \sum_{\ti{j}=2} z_0 z_j z_{\delta(j)} +
    \frac{1}{6} \sum_{\ti{j}=2, \ti{l}=2,\ti{k}=4} \CB(T_j * T_l,
    T_{\delta(k)}) z_j z_l z_{\delta(k)}.
$$
  In the special case where $\dim V_{2p} = 1$ for $p=0,\ldots,3$, we
  obtain $\phi_0(z_0,\ldots,z_3) = z_0 z_1 z_2 + \kappa z_1^3$ for
  $\kappa = \frac{1}{6} \CB(T_1*T_1,T_1)$. This is the case of, for
  instance, the central part of the cohomology of the quintic
  threefold in $\CP^4$.
\endexample

Next we consider deformations of a framed Frobenius module induced
from a potential. Let $R:=\C\{q_1,\ldots,q_r\}_0$ denote the ring of
convergent power series vanishing for $q_1 = \cdots = q_r = 0$ and
$R'$ be its image under the map induced by $q_j\mapsto e^{2\pi i z_j}$
for $1\leq j \leq r$.

\definition{Definition 4.8}
  Let $(V,\CB,e,*)$ be a framed Frobenius module of weight $k > 3$
  with classical potential $\phi_0$. A \jdef{quantum potential} on $V$
  is a function $\phi:V\rightarrow \C$ of the form $\phi = \phi_0 +
  \phi_\h$, where
  $$
    \phi_\h(z) \ :=\ \sum_{\ti{a}=2k-4} z_a
    \phi_h^a(z_1,\ldots,z_r) +\sum_{2 < \ti{a}< 2k-4,\,
    \ti{a}+\ti{b} = 2k-2} z_a z_b \phi_h^{a b}(z_1,\ldots,z_r), 
                                                               \tag{4.9}
  $$
  with $\phi_h^a, \phi_h^{a b}\in R'$ and such that the action 
  $$
    T_j \cdot_q T_a \ :=\ \sum_{\ti{c}=\ti{a}+2}
    \frac{\partial^3 \phi}{\partial z_j \partial z_a \partial
    z_{\delta(c)}} T_c\,,\quad \text{ with } q=(q_1,\ldots,q_r) \in
    \Delta^r\,. 
                                                                \tag{4.10}
  $$
  turns $(V,\CB,e,\cdot_q)$ into a graded $V_2$-Frobenius module for
  all $q$.
  
  For weight $k=3$ modules, the form of the potential is just $\phi =
  \phi_0 +\phi_\h$ where $\phi_\h\in R'$, and requiring that~(4.10)
  defines a graded $V_2$-module imposes no constraint on $\phi_\h$,
  contrary to the situation of higher weights.
\enddefinition

\remark{Remark} For $V$ of weight $1$ or $2$, the Frobenius
  module is determined by $\CB$ and $e$; hence no deformations in the
  sense of Definition~(4.8) are possible.
\endremark

\remark{Remark 4.11} 
  The form~(4.9) of the quantum
  potential is motivated by the Gromov-Witten potential, except that
  in the setup of Frobenius modules only a graded part of the
  potential is relevant.

  The condition for a quantum potential to induce a Frobenius module
  can be interpreted as a WDVV type of equation, only that, compared
  to the original setup in quantum cohomology, here we only see a
  graded part of the full system.
\endremark

\example{Example 4.12} In the weight $k=3$ case we have, for
  $\ti{j}=\ti{l}=2$,
  $$ 
    T_j \cdot_q T_l = T_j * T_l + (2\pi i)^3 \sum_{\ti{a}=4}
    \frac{\partial^3 \phi_\h(q_1,\ldots,q_r)}{\partial q_j \partial
    q_l \partial q_{\delta(a)}} q_j q_l q_{\delta(a)} T_a,
  $$ 
  while the action of $T_j$ on all the other graded parts of $V$ is
  the one given by $*$. Notice that, since $\phi_\h$ is assumed to be
  convergent at $q_1=\cdots=q_r=0$, we have $\cdot_0 = *$.
\endexample

As we stated at the beginning of this section, there is a
correspondence between deformations of a framed Frobenius module and a
certain type of maximally degenerating polarized variation of Hodge
structure. To make precise the notion of maximally degenerating
variation we reformulate Definition~(2.27) in the context of abstract
variations.

\definition{Definition 4.13}
  Given a polarized variation of Hodge structure of weight $k$ over
  $(\Delta^*)^r$ whose monodromy is unipotent, we say that $0\in
  \Delta^r$ is a \jdef{maximally unipotent boundary point} if
  \roster
  \item $\dim I^{k,k} = 1$, $\dim I^{k-1,k-1} = r$ and $\dim I^{k,k-1}
    = \dim I^{k-2,k}= 0$, as well as $I^{p,q} = 0$ for all $p,q<0$,
    where $I^{*,*}$ is the bigrading associated to the limiting mixed
    Hodge structure and,
  \item $\vspan_\C\{N_1(I^{k,k}),\ldots,N_r(I^{k,k})\} = I^{k-1,k-1}$,
    where $N_j$ are the monodromy logarithms of the variation.
  \endroster
\enddefinition

\proclaim{Theorem 4.14}
  There is a 1-1 correspondence between
  \roster
  \item"(a)" Deformations of a framed Frobenius module
    $(V,\CB,e,*)$ of weight $k$ arising from a quantum potential, and
  \item"(b)" Germs of polarized variations of pure Hodge structure of
    weight $k$ on $V$ degenerating at a maximally unipotent boundary
    point to a limiting mixed Hodge structure of Hodge-Tate type,
    split over $\R$, and together with a marked real point $e\in F^k$.
  \endroster
\endproclaim

The proof of this correspondence can be found in~\cite{CF2, Theorem
4.1}. The cases of weight $3$ and $4$ have been analyzed in~\cite{P}
and~\cite{CF1}. Below we will describe the main constructions that set
up the correspondence. The key technical step in proving the theorem
is the asymptotic description of the Hodge filtration of an admissible
variation presented in \S 3.

Let $\V\to (\Delta^*)^r$ be a variation of pure, polarized Hodge
structure with unipotent monodromy.  Then, by Corollary (3.20), the
germ of $\V$ at zero is determined by the following data:
\roster
\item The nilpotent orbit $(3.12)$;
\item The function $\psi(s) = e^{\Gamma(s)}\cdot F_{\infty}$.
\endroster

Moreover, to such a nilpotent orbit $(3.12)$, we can associate the
corresponding limiting mixed Hodge structure $(F_{\infty},W[-k])$,
which turns out to be {\it polarized} by the monodromy cone $(3.14)$
(cf\. \cite{CK}).  Furthermore, this polarized mixed Hodge structure
is equivalent to the nilpotent orbit~\cite{Sc}, \cite{CKS},
\cite{CF1,Theorem 2.3}.  Also,~\cite{CF1, Theorem 2.7}, on account of
the horizontality of the period map of $\V$, the function $\Gamma(s)$
can be recovered from its projection $\Gamma_{-1}$ onto the subspace
$\jgp_{-1} := \oplus_s\, \gg^{-1,s}$ of [cf\. Lemma $(3.17)$]
$$
     \gg_- := \qpsi = \bigoplus_{r<0,\, s}\, \gg^{r,s}
$$
and the monodromy logarithms $N_1,\dots, N_r$ of $\V$. 

The period mapping of the variation $\V$ is written as
$e^{X(s)}\cdot F_\infty$ for a (multivalued) map
$X:(\Delta^*)^r\rightarrow \gg_-$, as in (3.29)-(3.30). The
horizontality condition of the variation $\V$ can be written, in terms of
$X$, as
$$
    e^{-X(s)}\, d(e^{X(s)}) \ =\ dX_{-1} \quad\text{ for } X_{-1} \ =\
    \sum_{j=1}^r \frac{\log(s_j)}{2\pi i} N_j + \Gamma_{-1}(s),
                                                                \tag{4.15}
$$
and, in turn, this is equivalent to $dX_{-1}$ being a Higgs field.

\remark{Remark} 
  Theorem~(4.14) can now be refined to show that under the same
  correspondence, the nilpotent orbit of the variation of Hodge
  structure corresponds to the framed Frobenius module and that the
  function $\Gamma_{-1}$ corresponds to the quantum
  potential. Moreover, the transversality condition of the variation
  is equivalent to the graded part of the WDVV equation alluded to in
  Remark~(4.11).
\endremark

In what follows we will describe the two constructions that establish
the equivalence between degenerating variations of Hodge structure and
Frobenius modules.

\example{Construction 4.16 (Variation of Hodge structure from a
  Frobenius module)} This is the transcription of the well known
  $A$-model variation of Hodge structure into the language of families
  of Frobenius modules.  See \cite{CF2, \S 5} for more details as well
  as~\cite{CoK, Chapter 8} for the standard $A$-model variation.

  Consider a deformation of the framed Frobenius module $(V,\CB,e,*)$
  of weight $k$ generated by the potential $\phi = \phi_0 + \phi_\h$,
  with the framing $\{T_1,\ldots,T_r\}$; let $(q_1,\ldots,q_r)$ be the
  corresponding system of local coordinates on $(\Delta^*)^r$. Define
  \roster
  \item A free sheaf $\V := V \otimes \O_{(\Delta^*)^r}$.
  \item A Hodge filtration $\F^p := (\oplus_{a\geq p} V_{2(k-a)}) \otimes
    \O_{(\Delta^*)^r} \subseteq \V$. 
  \item The Dubrovin connection: for $T\in V$ and $j=1,\ldots r$,
    $$
      \nabla_{\jpd{q_j}} T \ :=\ \frac{1}{2\pi i q_j} T_j \cdot_q T.
    $$
    It follows from the symmetry condition (4.3) that $\nabla$ is
    flat. Also, it has simple poles at $q_j = 0$ where the residue
    $\Res_{q_j=0}(\nabla)$ is, up to a constant, the endomorphism of
    $(V,\CB,e,*)$ given by the action of $T_j$ on $V$. The logarithms of
    the monodromy of $\nabla$, $N_j$, are also given by that same action.
  \item A polarization $\CQ$ given by (4.5).
  \item A real form $\V_\R$ defined as follows. Let 

    $$
      \ti{\V} = V \otimes \O_{\Delta^r} \text{ and } 
      \ti{\nabla} = \nabla - \frac{1}{2\pi i}\sum_{j=1}^r N_j\otimes
      \frac{dq_j}{q_j}.
    $$ 
    Then $\ti{\nabla}$ is a flat connection on $\ti{\V}$; for $v\in V$
    we define $\ti{\sigma}_v$ to be the $\ti{\nabla}$-flat section
    such that $ \ti{\sigma}_v(0) = v$. Then $\V_\R \subseteq \V$ is the
    local system generated by the sections $\exp(-\frac{1}{2\pi
    i}\sum_j\log(q_j)N_j)\ti{\sigma}_v$ for all $v\in V_\R$.
  \endroster
\endexample

\proclaim{Theorem 4.17}
  The tuple $(\V,\nabla,\F,\V_\R,\CQ)$ constructed above is a
  polarized variation of Hodge structure of weight $k$ having a
  maximally unipotent boundary point at $0\in \Delta^r$ where the
  mixed Hodge structure is of Hodge--Tate type split over $\R$.
\endproclaim

A proof of Theorem~(4.17) can be found in~\cite{CF2, Theorem 5.10}.
See~\cite{CoK, Theorem 8.5.11} for a version of this result in the original
setup of quantum cohomology of Calabi-Yau manifolds.

\remark{Remark}
  The construction described above does not explicitly use the fact
  that the deformation of the Frobenius module comes from a
  potential. Nevertheless, the quoted proof does depend on that fact. 
\endremark

\example{Construction 4.18 (Frobenius module from a variation of Hodge
  structure)} The starting point is a polarized variation of Hodge
  structure of weight $k$ degenerating to a maximally unipotent
  boundary point where the limiting mixed Hodge structure is of
  Hodge--Tate type, split over $\R$, together with a real element $e\in
  F^k$.

  The construction of a Frobenius module uses only the information of
  the variation near the degeneration, so that we can restrict all
  considerations to the case of a variation over $(\Delta^*)^r$ having
  the origin as a maximal degeneration.

  Define:
  \roster
  \item A graded vector space $V:= \oplus_{p=0}^k V_{2p}$, with
    $V_{2p}$ defined by~(4.4).
  \item A symmetric nondegenerate pairing $\CB$ on $V$ by~(4.5).
  \item A real structure $V_\R$ on $V$ defined by the real structure
    carried, by hypothesis, by all $I^{p,p}$.
  \item A framing $\{T_1,\ldots,T_r\}$ of $V$, where $T_j:= N_j(e) \in
    V_2$ for $j=1,\ldots,r$.
  \item An action $*$ of $V_2$ on $V$ by $T_j * T:= N_j(T)$ for
    $j=1,\ldots,r$. 
  \item A ``unit'' $e \in F^k = I^{k,k} = V_0$.
  \endroster

  Then, $(V,\CB,*,e)$ is a $V_2$ Frobenius module framed by
  $\{T_1,\ldots,T_r\}$. This follows from the commutativity of the
  monodromy logarithms, the fact that these are infinitesimal
  isometries of $Q$, and that they polarize the limiting mixed Hodge
  structure.

  To define a deformation of $(V,\CB,*,e)$, we need special
  coordinates $(q_1,\ldots q_r)$ on $\Delta^r$, known as
  \jdef{canonical coordinates}. Such a coordinate system is
  characterized by the fact that the function $\Gamma_{-1}$ is
  normalized with respect to these coordinates to satisfy
  $\Gamma_{-1}(I^{1,1})=0$. Canonical coordinates exist because of the
  maximal unipotency condition on the variation (see~\cite{CF1, \S
  3}). Finally, for $j=1,\ldots,r$ and $T\in V$, define
  $$
    T_j\cdot_q T := \left(\frac{\partial X_{-1}}{\partial z_j}\right)_q (T)
                                                     \tag{4.19}
  $$
  where the right hand side is the evaluation of the endomorphism of
  $V$, $\frac{\partial X_{-1}}{\partial z_j}$, on the vector
  $T$. The fact that this action defines a family of $\Sym V_2$
  modules is equivalent to~(4.15). The quantum potential $\phi_\h$
  that generates the deformation is:
  $$
    \split
    \phi_h^{a b}(q) &:= \frac{1}{2}\, \CB(\Gamma_{-1}(T_a),T_b)
    \text{ for } 2<\ti{a}< 2k-4 \text{ and } \ti{a}+\ti{b}=2k-2\\  
    \phi_h^a(q) &:= \CB(-\Gamma_{-2}(T_a),T_0) \text{ for }
    \ti{a}=2k-4\\
    \phi_\h &:= \sum_{\ti{a}\ =\ 2k-4} z_a \phi_h^a +
    \sum_{2\ <\ \ti{a}\ <\ 2k-4,\, \ti{a}+\ti{b}\ =\ 2k-2}
    z_a z_b \phi_h^{a b}\\
    \phi  &:= \phi_0 + \phi_\h,
    \endsplit
                                                \tag{4.20}
  $$
  where $\Gamma_{-2}$ is determined by $\Gamma_{-1}$ and the monodromy
  logarithms. 
\endexample

\proclaim{Theorem 4.21}
  $(V,\CB,e,*)$ as constructed above is a graded Frobenius module
  framed by $\{T_1,\ldots,T_r\}$. Moreover, the $V_2$ action defined
  by~(4.19) is a deformation of $(V,\CB,e,*)$ generated by
  the potential~(4.20), whose quantum part is uniquely determined by the
  function $\Gamma_{-1}$ of the variation.
\endproclaim

\remark{Remark}
  In the weight $k=1$ and $k=2$ cases, the use of canonical
  coordinates implies the normalization $\Gamma = 0$. This choice
  parallels the fact that there are no deformations of Frobenius
  modules for the same weights. 
\endremark


\head 5.\quad Frobenius Algebras
\endhead

Following Dubrovin, a \jdef{Frobenius algebra} over the finite
dimensional $\C$-vector space $V$ consists of a tuple $(V,\CB, \apc,
e)$, where $\apc$ defines a commutative and associative $\C$-algebra
structure on $V$ with unit $e$, and $\CB$ is a non-degenerate
symmetric bilinear form on $V$ such that
$$
  \CB(v_1\apc v_2, v_3) \ =\ \CB(v_2, v_1 \apc v_3) \text{ for all } v_1,
  v_2, v_3 \in V.                                       \tag{5.1}
$$ 
Additionally, $(V,\CB, \apc, e)$ is a \jdef{graded Frobenius algebra
of weight $k$} if $V$ decomposes as $V=\oplus_{p=0}^k V_{2p}$, $V_0 =
\langle e\rangle$, the product $\apc$ is graded and $\CB$ pairs
$V_{2p}$ with $V_{2(k-p)}$ for all $p$.

If $(V,\CB,\apc,e)$ is a graded Frobenius algebra of weight $k$ we
can, by restricting the action of $\apc$, define a graded $V_2$-module
$(V,\CB,*,e)$.  In general, this process can not be reversed since
part of the product structure is lost. Nevertheless, in some
cases a full Frobenius algebra can be constructed from a graded
$V_2$-module. We will illustrate this phenomenon in the following
examples.

\example{Example 5.2} Let $(V,\CB,*,e)$ be a graded $V_2$ Frobenius
  module of degree $k=3$. Then, for homogeneous elements $v_a$, $v_b$
  in $V_{\ti{a}}$ and $V_{\ti{b}}$ we define a multiplication $\apc$
  by: $v_a \apc e = e \apc v_a = v_a$; for $\ti{a} = 2$ ($\ti{b}$
  arbitrary), $v_b \apc v_a = v_a \apc v_b = v_a * v_b$ and, $v_a \apc
  v_b = 0$ if $\ti{a}>2$ and $\ti{b} >2$. We notice that $\apc$ is
  graded and commutative.  Associativity can be seen as follows: if
  any vector in a triple product is (a multiple of) the identity,
  associativity is immediate. If all vectors are in $V_2$,
  associativity follows from~(4.2).  Finally, if none of these
  conditions hold, all triple products vanish by the grading
  properties. Finally, a short analysis shows that $\CB(v_a \apc v_b,
  v_c) = \CB(v_b, v_a \apc v_c)$.

  Similar arguments show that the same result applies in the weight
  $k=4$ and $k=5$ cases. The cases of weight $k=1$ and $k=2$ are
  trivial since the product structure is determined by $\CB$ and $e$.
\endexample

\example{Example 5.3} 
  Let $(V,\CB,*,e)$ be a graded $V_2$ Frobenius module of degree $k$.
  We say that \jdef{$V$ is generated by $V_2$} if the linear map $\Sym
  V_2 \rightarrow V$ defined by $P\mapsto P*e$ is onto. In this
  case, the algebra structure of $\Sym V_2$ can be transferred to $V$.
  Explicitly, for homogeneous vectors $v_a = P_a * e$ and $v_b = P_b *
  e$ we have $v_a \apc v_b = (P_a \cdot P_b) * e$. Clearly, $P_a$ is
  homogeneous, so that $\apc$ defines a graded multiplication that is
  commutative, associative and with unit $e$.  Condition~(5.1) follows
  immediately from the iteration of~(4.2).

  If $(V,\CB,\cdot_q,e)$ is a (continuous) family of Frobenius modules
  such that $\cdot_0 = *$ generates $V$ in $V_2$, then for all $q$ near
  $0$ the corresponding Frobenius module also generates and we obtain
  a family of graded Frobenius algebras defined over an open
  neighborhood of $q=0$ of the same parameter space.
\endexample

In quantum cohomology the standard constructions produce the Frobenius
manifold structure on $H^*(X,\C)$ whose product is the \jdef{big
quantum product}. Then, considering $i:H^2(X,\C)\subset H^*(X,\C)$,
the product obtained by pulling back to $i^*(T H^*(X,\C)) \simeq H^2
\times H^*$ the big quantum product is known as the \jdef{small
quantum product}.

Starting from a family of graded $V_2$-Frobenius modules we have seen
that in some cases it is possible to construct a family of Frobenius
algebras over an open subset of the same base. In Section~4 we took
the parameter space of such a family to be $\Delta^r$, with $r = \dim
V_2$. We can also consider the covering of the polydisk by the product
of upper half planes, $U^r$, given by $z\mapsto q = \exp(2\pi i z)$;
in this case we pull back our constructions and obtain a family of
modules (or algebras) over $U^r$, that we view as contained in $V_2$
via the framing. In this case, the family is defined near the point at
infinity.

One can pose then the following question: If $V$ is generated by $V_2$
under $*$, the construction described in the previous paragraph
resembles the small quantum product with a family of Frobenius
algebras defined over (an open set contained in) $V_2$. {\it Is it
possible to ``unfold'' this structure to obtain a full Frobenius
manifold analogous to the big quantum product?} We will show that the
answer to this question is positive in some cases.

Before going into the details, recall that a \jdef{Frobenius
manifold}~\cite{Du} is a complex manifold $M$ of dimension at least
$1$ with a commutative and associative multiplication $\apq$ on the
holomorphic tangent bundle $TM$, a unit field $e$ and a symmetric,
nondegenerate bilinear pairing $g$ on $TM$ such that the Levi-Civita
connection of $g$, $\nabla^g$, is flat, the unit field $e$ is
$\nabla^g$-flat and for all vector fields $X$, $Y$, $Z$ on $M$ the
following conditions hold:
$$
  g(X\apq Y, Z) \ =\ g(Y, X \apq Z)                     \tag{5.4}
$$
and
$$
  \nabla^g_X(Y\apq Z) - \nabla^g_Y(X\apq Z) + X\apq \nabla^g_Y Z - 
  Y \apq \nabla^g_X Z - [X,Y] \apq Z \ =\ 0. \tag{5.5}
$$
\remark{Remark}
  Frobenius manifolds usually carry a vector field known as the
  ``Euler field''. We will not describe this field in our
  construction below, but such object can readily be constructed in terms
  of the grading that is part of a Frobenius module.
\endremark

\example{Construction 5.6}
  Let $(V,\CB,*,e)$ be a framed Frobenius module of weight $k$ that is
  deformed via the potential $\phi = \phi_0 + \phi_\h$ and assume that
  the module structure of $*$ can be extended to an algebra that we
  denote by $\apc$. As we mentioned above the family can be lifted to a
  family of modules on $U^r\subseteq V_2$. An adapted basis $T_0,\ldots
  T_m$ of $V$ provides coordinates $z_0,\ldots z_m$ and the algebra
  structure $\apc$ corresponds to $\Im(z_j)=\infty$ for $j=1,\ldots,
  r$.

  The algebra structure $\apc$ can be encoded into a classical
  potential $\hat{\phi}_0$, and we define a multiplication on the
  tangent bundle $TV$ by 
  $$
    T_a \apq_z T_b \ :=\ \sum_{c} \frac{\partial^3 (\hat{\phi}_0 +
    \phi_\h)}{\partial z_a \partial z_b \partial z_{\delta(c)}} T_c \ =\ 
     T_a \apc T_b + \sum_c \frac{\partial^3
    \phi_\h}{\partial z_a \partial z_b \partial z_{\delta(c)}} T_c,
                                                \tag{5.7}
  $$
  for all $\sum_a z_a T_a \in V$. We remark that this definition does
  not, in general, define a {\it graded} product.

  We also consider the (constant) metric on $TV$ induced by $\CB$ and
  its Levi-Civita connection $\nabla$ characterized by $\nabla T= 0$
  for all $T\in V$ regarded as a constant vector field. In particular,
  the unit $e \in V_0$ extends to a flat vector field that we keep
  denoting by $e$.

  To see if the previous data makes $TV$ into a Frobenius manifold we
  have to check the following
  \roster
  \item Compatibility between $\apq$ and $\CB$~(5.4): 
      $$
        \split
        \CB(T_a \apq_z T_b,T_c) &= \CB(T_a \apc T_b + \sum_{d}
        \frac{\partial^3 \phi_\h}{\partial z_a \partial z_b \partial
        z_{\delta(d)}} T_d, T_c) \\ &= \CB(T_a \apc T_b, 
        T_c) + \sum_{d}\frac{\partial^3 \phi_\h}{\partial z_a z_b
        z_{\delta(d)}} \CB(T_d,T_c) \\&= 
        \CB(T_a \apc T_b, T_c) + \frac{\partial^3 \phi_\h}{\partial
        z_a \partial z_b \partial z_c},
        \endsplit
      $$
    from where~(5.4) follows by applying~(5.1) to $\apc$
    and by the symmetry of the second summand.
  \item The potentiality condition~(5.5): 
      $$
      \split
        \nabla_{\jpd{z_a}}(T_b \apq_z T_c) &= \nabla_{\jpd{z_a}}(T_b
        \apc T_c + \sum_d \frac{\partial^3 \phi_\h}{\partial z_b
      \partial z_c \partial z_{\delta(d)}} T_d) \\ &=
        \nabla_{\jpd{z_a}}(T_b \apc T_c)  +
        \nabla_{\jpd{z_a}}(\sum_d \frac{\partial^3 \phi_\h}{\partial
        z_b \partial z_c \partial z_{\delta(d)}}
        T_d) \\&= \sum_d \frac{\partial^4 \phi_\h}{\partial z_a
        \partial z_b \partial z_c 
        \partial z_{\delta(d)}} T_d.
      \endsplit
    $$ 
    By the linearity of~(5.5) it suffices to check it for $X=T_a$,
    $Y=T_b$ and $Z=T_c$. The result then follows from the symmetry of
    the last expression.
  \item Commutativity, associativity and unit: commutativity is clear
    from formula~(5.7). That $e$ is a unit follows from
    $e$ being a unit for $\apc$ and that there is no
    dependence on $z_0$ in $\phi_\h$. The associativity of $\apq_z$,
    will be discussed below.
  \endroster

  We will only check the associativity in the cases of weight $3$, $4$
  and $5$. The fact that $\apq_z$ is not graded makes the computations
  hard, and this is only partially eased by the weaker property
  $V_{2p} \apq_z V_{2q} \subseteq \oplus_{a\geq (p+q)} V_{2a}$.

  {\sl Weight $k=3$}. In this case $\apq_z$ is graded and can be
  computed from
  $$
    T_a \apq_z T_b = 
      \left\{\aligned
       & T_a\cdot_{z'} T_b \text{ if } \ti{a} = \ti{b} = 2,\\
       & T_a \apc T_b \text{ otherwise},
      \endaligned
      \right.
  $$
  where $z' = \pi_2(z)$ is the projection of $\sum z_a T_a$ on $V_2$. 
  Associativity follows immediately.
  
  {\sl Weight $k=4$}. The more involved of the triple products is that of
  $T_a \apq_z (T_b \apq_z T_c)$ with $\ti{a}=\ti{b} = \ti{c} =2$. We
  proceed as follows. Remember that the quantum potential for weight
  $4$ is $\phi_\h = \sum_{\ti{a}=4} z_a \phi^a_h(z_1,\ldots,z_r)$.
  Then:
  $$
    \split
      T_a \apq_z (T_b \apq_z T_c) &= T_a \apq_z (T_b * T_c +
      \sum_{\ti{d}=4} \frac{\partial^3 \phi_\h}{\partial z_b \partial
    z_c \partial z_{\delta(d)}} T_d +
      \sum_{\ti{d}=6} \frac{\partial^3 \phi_\h}{\partial z_b \partial
    z_c \partial z_{\delta(d)}} T_d) \\&= T_a
      \apq_z (T_b \cdot_{z'} T_c + \sum_{\ti{d}=6}
      \frac{\partial^3 \phi_\h}{\partial z_b \partial z_c \partial
    z_{\delta(d)}} T_d) \\&= T_a \apq_z (T_b 
      \cdot_{z'} T_c) + \sum_{\ti{d}=6}
      \frac{\partial^3 \phi_\h}{\partial z_b \partial z_c \partial
    z_{\delta(d)}} T_a * T_d \\&= T_a\cdot_{z'} (T_b 
      \cdot_{z'} T_c) + \frac{\partial^3 \phi_\h}{\partial z_b
    \partial z_c \partial z_a} T_{\delta(0)}.
    \endsplit
  $$
  where, as before, $z'=\pi_2(z)$ denotes the projection of $z$ to
  $V_2$. Then, using the commutativity of $\apq_z$ we have
  $$ 
    (T_a \apq_z T_b) \apq_z T_c = T_c \apq_z (T_a \apq_z T_b) = T_c
    \cdot_{z'} (T_a \cdot_{z'} T_b) + \frac{\partial^3
    \phi_\h}{\partial z_a \partial z_b \partial z_c} T_{\delta(0)},
  $$
  and, since $\cdot_{z'}$ defines a family of Frobenius modules,
  relation~(4.3) implies the associativity of $\apq_z$.
  
  The {\it weight $k=5$} case is similar to the previous one, with more
  painful computations.
\endexample

\remark{Remark} 
  All together, the construction described in~(5.6) provides
  unfoldings of the original algebra structure defined on $V_2\times
  V$ to a full Frobenius manifold on $TV$ for weights $k=3, 4, 5$.  We
  note that these results do not follow from the application
  of~\cite{HM, Theorem 4.5} because in our case we are working with
  germs at infinity or, in terms of families defined over $\Delta^r$,
  our Higgs fields have logarithmic singularities. Perhaps,
  Construction~(5.6) should be regarded as evidence for an extended
  version of~\cite{HM, Theorem 4.5}.

  It seems unlikely that associativity ---hence this explicit
  construction--- can be extended to weights $k\geq 6$. The reason for
  that is that the homogeneity properties of a potential defining the
  product on a Frobenius manifold imply that, for weight $k\leq 5$ all
  potentials have the form~(4.9), whereas for $k\geq 6$ new homogeneous
  terms can be present.
\endremark

\vskip .5cm

Another connection between families of Frobenius modules generated in
$V_2$ and Frobenius manifolds can be established as follows. 

\proclaim{Theorem 5.8} 
  Consider a deformation of the framed Frobenius module $(V,\CB,*,e)$
  induced by a quantum potential defined over a neighborhood of $0\in
  \Delta^r$ $(r = \dim V_2)$ and such that $V$ is generated by $V_2$
  under $*$.  Then, for each $\hat{s}\in (\Delta^*)^r$ near $0$ the
  family of Frobenius modules generates a germ of a family of algebras
  at $\hat{s}$.  Moreover, this germ can be unfolded to a germ of a
  Frobenius manifold.
\endproclaim

\demo{Proof} 
  By Theorem~(4.14) the deformation of $(V,\CB,*,e)$ generates a
  polarized variation of Hodge structure $(\V,\nabla,\F,\V_\R,\CQ)$
  over a neighborhood of $0\in \Delta^r$, with a marked element $e \in
  F^k$. The limiting mixed Hodge structure of this variation is given
  by~(4.4). As was observed in~(3.13), this variation defines an
  admissible variation of graded-polarized mixed Hodge structure with
  unipotent monodromy.

  As in~(3.16) let $\Psi_p = \oplus_{a\leq p,\,b} I^{a,b} =
  \oplus_{a\leq p} V_{2(k-a)}$. Then $\Psi_p$ is $\gg_-$-invariant and
  opposite to $F_\infty$. Moreover, by Theorem~(3.28) $\Psi_p$ extends
  to a $\nabla$-flat increasing filtration $\underline{\Psi}_p$ of
  $\V$ that pairs with $\F^p$ to define an unpolarized $\C$VHS in a
  neighborhood of the origin. In particular, by~(1.3), $\F^p$ and
  $\underline{\Psi}_p$ are opposite.

  By parallel transporting the bundles and additional structure to a
  fixed fiber of $\V$ we write
  $$
    F^p(s) = e^{X(s)}\cdot F^p_\infty \text{ with } F^p_\infty =
    \oplus_{a\geq p} V_{2(k-p)}.
  $$
  Recalling that $\Psi_p$ is $\gg_-$-invariant we can find the
  $C^\infty$ decomposition~(1.4) of the complex variation as follows:
  $$
    \split
    \U^p(s) &= F^p(s) \cap \Psi_p(s) = (e^{X(s)}\cdot F^p_\infty) \cap
    \Psi_p = e^{X(s)}\cdot (F^p_\infty \cap \Psi_p) \\ &= e^{X(s)}\cdot
    V_{2(k-p)}.
    \endsplit 
  $$
  Then, using~(1.5), the Higgs field $\theta$ of the $\C$VHS is 
  $$
    \theta(s) \ =\ \exp(X(s))\, dX_{-1}\, \exp(-X(s)).           \tag{5.9}
  $$
  
  We notice that since $X(s)$ is an infinitesimal automorphism of $Q$,
  $$
    \CQ(\F^p, \F^{k-p+1}) = Q(F^p_\infty,F^{k-p+1}_\infty) = 0
  $$
  and, by the $\gg_-$-invariance of $\Psi_p$, 
  $$
    \CQ(\underline{\Psi}_p(s),\underline{\Psi}_{k-p-1}(s)) =
    Q(\Psi_p,\Psi_{k-p-1}) = 0.
  $$

  From~(5.9) it follows that, for a fixed $\hat{s}$, the $\C$-span of the
  iterated application of $\theta_Y$ for $Y\in T_{\hat{s}} (\Delta^*)^r$ is
  conjugate to that of $(dX_{-1})_Y$. But the span for these last
  operators as $\hat{s}\rightarrow 0$ is the span of the iterated action of
  the monodromy logarithms which act as $V_2$ acts on $V$ via $*$ as
  was remarked in (3) of Construction~(4.16). Then, under the
  hypothesis that $V$ is generated by $V_2$ under $*$, we conclude
  that the same condition holds for the Higgs field $\theta(\hat{s})$, for
  $\hat{s}$ near $0$.

  Next we use the machinery of~\cite{HM, \S 5} to construct a
  Frobenius manifold. We claim that 
  $(((\Delta^*)^r,\hat{s}),\V,\nabla,\F)$ is a germ of an
  $H^2$-generated variation of filtrations of weight $k$, in the
  language of~\cite{HM, Definition 5.3}. The only thing that remains
  to be proved is that the Higgs field $C:\F^p/ \F^{p+1} \rightarrow
  \F^{p-1}/\F^p$ induced by $\nabla$ satisfies the generation
  condition. But, under the isomorphism introduced by parallel
  transport to a fixed fiber of $\V$ and the canonical isomorphism
  $F^p(s)/ F^{p+1}(s) \simeq U^p(s)$ we see that $C$ is, in fact,
  isomorphic to $\theta$ and, since this last field satisfies the
  generation condition, so does $C$, and the claim is proved.

  But then, the $H^2$-generated variation of filtrations of weight $k$
  defined above, together with the pairing $\CQ$, the opposite
  filtration $\underline{\Psi}_p$ and the ``fixed generator'' $e$
  satisfy the conditions of~\cite{HM, Theorem 5.6}. Hence, we conclude
  that there is a unique unfolding of this structure to a germ of a
  Frobenius manifold.
\enddemo

\remark{Remark} 
  The filtration $\Psi$ constructed in Theorem $(5.8)$ is related to
  the relative weight filtration $\rel W$ by the rule $\Psi_p = \rel
  W_{2p}$ since the associated limiting mixed Hodge structure is
  Hodge--Tate.
\endremark

In view of Theorem~(4.14) there is the following immediate corollary
to Theorem~(5.8).

\proclaim{Corollary 5.10}
  Let $\V\rightarrow (\Delta^*)^r$ be a variation of Hodge structure
  of weight $k$ which satisfies the hypothesis of Theorem $(4.14)$,
  and assume that the iterated action of the monodromy cone of $\V$ on
  $e\in F^k_\infty \subseteq F^0_\infty$ spans $F^0_\infty$. Then, for
  each $\hat{s}\in\Delta^{*r}$ sufficiently close to the origin, there
  exists a corresponding germ of a Frobenius manifold $M_{\hat{s}}$,
  which is completely determined by the asymptotic behavior of $\V$.
\endproclaim

\remark{Remark}
  Because of~(5.9) we see that there is a simple
  connection between the Higgs field $\theta$ and the one that appears
  in the correspondence described in
  Section~4, namely, $dX_{-1}$. In fact, 
  $$
    \theta \ =\ \exp(X)\, dX_{-1}\, \exp(-X) \ \equiv\ dX_{-1} 
    \mod F^{-2} \gg.
  $$
  Alternatively, $\theta$ and $dX_{-1}$ define isomorphic Higgs
  bundles. 
\endremark


\enddocument